\documentclass[a4paper,12pt]{article}
\usepackage[T1]{fontenc}
\usepackage[cp1250]{inputenc}
\usepackage[pdftex]{graphicx}
\usepackage{amsmath,amstext,amssymb,amsopn,amsthm,amsfonts,mathrsfs,extraipa,dsfont,bbm}

\newcommand{\I}{\mathcal{I}}

\newcommand{\Ip}{\mathcal{I}^+}

\newcommand{\In}{\leqslant_{\tiny\mathcal{I}}}

\newcommand{\D}{\mathcal{D}}

\begin{document}

\title{{\bf Remarks on quotient algebras}}
\author{Ryszard Frankiewicz \and S{\l}awomir Szczepaniak}
\date{}
\maketitle

\begin{abstract}
The structure of quotient Boolean algebras in terms of cardinal invariants is investigated. Some results of Gitik and Shelah regarding atomless ideals are reproved and proofs are significantly simplified.
\end{abstract}

{\bf 1. Introduction}

\vspace{0.5cm}

\indent While studying the structure of Boolean algebras one faces to the problem of classifying them up to (boolean) isomorphism. There are several ways of distinguishing between Boolean algebras, e.g. one can try to show than they produce (elementarily) nonequivalent Boolean extensions. This approach was used by Gitik and Shelah in \cite{GS}. Although elegant their method is rather indirect and hide a very reason why the given Boolean algebras are not isomorphic or equivalent (c.f. below). It is our opinion that more traditional way - by invariants - is better suited for the mentioned problem. The presented paper follows this approach. As a result we obtain not only much simpler proofs of Gitik-Shelah's results but also (as we believe) more illuminating. The main contribution are elementary (classical) proofs of {\bf\em Corollary 3} and {\bf\em Corollary 7}.

\indent For undefined notions we refer to \cite{J} and fix those most frequent in the paper or slightly different from the mentioned source. By an {\em ideal} on a given set $X$ we shall understand the family $\I\subseteq P(X)$ such that
\begin{itemize}
  \item if $A_0,A_1\in\I$ and $A\subseteq A_0\cup A_1$ then $A\in\I$;
  \item $X\notin\I$ and $\{x\}\in\I$ for each $x\in X$.
\end{itemize}
We denote by $\Ip:=\{A\subseteq X:A\notin\I\}$ a family of $\I$-{\em positive sets} and by $\I\upharpoonright Y:=\{A\subseteq X:A\cap Y\in\I\}$ an ideal called {\em restriction} of $\I$ to the set $Y$. An ideal $\I$ is $\kappa$-{\em complete} if is closed under sums of length less than $\kappa$. The {\em additivity number} of $\I$ is defined as
$$\mbox{add}(\I):=\min\left\{\kappa:\I\,\,\mbox{is not}\,\,\kappa\mbox{-complete}\right\}.$$ The additivity number is that it is always regular cardinal.

\indent We are mostly focused on the quotient Boolean algebra $P(X)/\I$ with the ordering
$$[A]\leqslant[B]\quad\mbox{iff}\quad A\In B\quad\mbox{iff}\quad A\setminus B\in\I.$$
There are several properties and characteristics of $\I$ formulated in terms of boolean structure of $P(X)/\I$. We say that an ideal $\I$ on $X$ is {\em atomless} if $P(X)/\I$ is an atomless Boolean algebras. In other words for no $A\in\Ip$ a restriction $\I\upharpoonright A$ is a {\em prime} ideal, i.e. $P(X)=(\I\upharpoonright A)\cup(\I\upharpoonright A)^c$. The paper deals only with atomless ideals as they are the most important, they produce nontrivial boolean extensions. In the latter case, remind, if two Boolean algebras have isomorphic dense subsets, we call them {\em equivalent} (as forcing notions) or {\em co-absolute}. The other idea connecting an ideal and its quotient Boolean algebra is the notion of saturation. We say that a Boolean algebra $\mathbbm{B}$ is {\em $\kappa$-saturated} if there is no antichain of size $\kappa$ in $\mathbbm{B}$ and define the {\em saturation number} of $\mathbbm{B}$ as
$$\mbox{sat}(\mathbbm{B}):=\min\left\{\kappa:\mathbbm{B}\,\,\mbox{is}\,\,\kappa\mbox{-saturated}\right\}.$$
Then an ideal $\I\subseteq P(X)$ is {\em $\kappa$-saturated} if $P(X)/\I$ is $\kappa$-saturated. Recall a Boolean algebra $\mathbbm{B}$ is called {\em $\kappa$-complete} if $\sum M$ exists for all subsets $M\subseteq\mathbbm{B}$ of size less than $\kappa$. If $\I$ is $\kappa$-complete then also $P(X)/\I$ is $\kappa$-complete. A subalgebra $\mathbbm{A}\subseteq\mathbbm{B}$ is a {\em regular subalgebra} of $\mathbbm{B}$ if for each $M\subseteq\mathbbm{A}$ such that $\sum_{\mathbbm{A}} M$ exists, also $\sum_{\mathbbm{B}} M$ exists and $\sum_{\mathbbm{A}} M=\sum_{\mathbbm{B}} M$. Thus, in this case, we can omit subscripts in $\sum$. We use a term {\em partition} for maximal antichain in $\mathbbm{B}$. A family $\D\subseteq\mathbbm{B}^+:=\mathbbm{B}\setminus\{\mathds{O}\}$ is {\em dense} in $\mathbbm{B}$ if for all $b\in\mathbbm{B}^+$ one can find $d\in\D$ such that $d\leqslant_{\tiny\mathbbm{B}}b$. We define the {\em density number} of $\mathbbm{B}$ (see \cite{M}) and {\em hereditary density number} as follows
$$\pi(\mathbbm{B}):=\min\left\{|\D|:\D\,\,\mbox{is dense in}\,\,\mathbbm{B}\right\},$$
$$h\pi(\mathbbm{B}):=\min\left\{\pi\left(\{b\in\mathbbm{B}:b\leqslant a\}\right):a\in\mathbbm{B}^+\right\}.$$
In view of an isomorphism $P(X)/(\I\upharpoonright A)$ and $P(A)/(P(A)\cap\I)$ we have
$$h\pi(P(X)/\I)=\min\left\{\pi\left(P(X)/(\I\upharpoonright A):A\in\Ip\right)\right\}.$$
The density of an ideal can be seen as the strong notion of saturation. For if there is an antichain of size $\kappa$ then a density has to be at least $\kappa$. An unpleasant feature of the density number is lack of dependence between density numbers of a Boolean algebra and its subalgebras. 

As our proofs are inspired by {\em Base Matrix Lemma} (\cite{BS},1.12) we need some notations concerning distributivity and trees. We say that a partition $\mathcal{P}_1$ {\em (almost) refines} a partition $\mathcal{P}_0$ if for all $a\in\mathcal{P}_1$ there is (finitely many) exactly one $b\in\mathcal{P}_0$ such that $a\cdot b\neq\mathds{O}$. A Boolean algebra $\mathbbm{B}$ is called (weakly) $\kappa$-distributive if for every family $\left(\mathcal{P}_\alpha\right)_{\alpha<\kappa}\subseteq P\left(\mathbbm{B}^+\right)$ of partitions there is a partition (almost) refining all of them. Define {\em (weakly) distributivity number} as
$$(w)h(\mathbbm{B}):=\min\left\{\kappa:\mathbbm{B}\,\,\mbox{is not (weakly)}\,\,\kappa\mbox{-distributive}\right\}$$
and observe that always $h(\mathbbm{B})\leqslant wh(\mathbbm{B})$. We abuse notation and write $\mbox{sat}(\I), \pi(\I), h(\I), wh(\I)$ whenever $\mathbbm{B}=P(X)/\I$. The last important feature of the above cardinal invariants is their forcing-invariance. This means that if $\theta\in\{\mbox{sat}, \pi, h, wh\}$ and a Boolean algebra $\mathbbm{A}$ is equivalent to $\mathbbm{B}$ then $\theta(\mathbbm{A})=\theta(\mathbbm{B})$.

By a tree $\mathcal{T}$ we shall understand a set of sequences of ordinals with the following ordering ('$\hat{}$' denotes concatenation)
$$s\leqslant_{\mathcal{T}}t\quad\mbox{iff}\quad s=t\,\hat{}\,w\,\,\,\mbox{for some sequence of ordinals}\,\,w$$
with the largest element being $\emptyset$ (empty sequences). An element $t\in\mathcal{T}$ is called {\em splitting node} (we say that $t$ splits) if in $\mathcal{T}$ there are two $\leqslant_{\mathcal{T}}$-incomparable extensions of $t$. We write $\left[\mathcal{T}\right]$ for a set of branches of $\mathcal{T}$ and $\mathcal{T}\upharpoonright\alpha$ for a subtree of $\mathcal{T}$ such that $$\mathcal{T}\upharpoonright\alpha=\bigcup_{\beta<\alpha}Lev_\beta(\mathcal{T})$$
where $Lev_\beta(\mathcal{T})$ stands for $\beta^{\mbox{\tiny th}}$ level of $\mathcal{T}$. The heights of a tree and its nodes is defined as usual, see e.g. \cite{J}. A family $\{b_t\in\mathbbm{B}:t\in\mathcal{T}\}$ indexed
by nodes of a tree $\mathcal{T}$ is called {\em Cantor (sub)scheme} if $b_{t_1}\leqslant_{\mathbbm{B}} b_{t_0}$ for $t_1\geqslant t_0$ and $b_{t_2}\cdot b_{t_3}=\mathds{O}$ for incomparable nodes $t_2, t_3$ of the tree $\mathcal{T}$.
%\indent In the paper we compare Boolean algebra $P(X)/\I$ with $C(\lambda)$ and measure algebra. The Boolean algebra $C(\lambda)$ is called {\em Cohen algebra} and is equivalent to forcing adding $\lambda$ many Cohen reals, i.e. it includes as a dense subset a poset of finite binary functions with domain being subset of $\lambda$ and ordering being a reverse inclusion. A boolean algebra $\mathbbm{B}$ is called {\em measure algebra} if it carries a measure, i.e. non-decreasing, real-valued, non-negative, $\sigma$-additive, normalized ($\mu(\mathbbm{1})=1$) function $\mu$ vanishing only at $\mathds{O}$. A measure $\mu$ is atomless iff $\I_\mu:=\{a\in\mathbbm{B}:\mu(a)=0\}$ is an atomless ideal.

\medskip

\noindent{\bf 2. Results}

\medskip

In this section we prove several assertions on structure of the atomless Boolean algebra $P(X)/\I$ by the means of cardinal invariants. Most of them appeared in \cite{GS} as results of their generic ultrapower technique. Our proofs are much simpler and in the spirit of 1.12 of \cite{BS}.

\medskip

\indent {\bf\em Theorem 1} If $\I$ is an atomless ideal on $X$ then $h\pi(\I)\geqslant\mbox{add}(\I)$.

\medskip

\indent  {\bf\em Proof:}

\medskip
Let $A\in\Ip$ be arbitrary and denote $\kappa:=\mbox{add}(\I)$. If $\mbox{sat}(\I\upharpoonright A)>\kappa$ then evidently $\pi(\I\upharpoonright A)\geqslant\kappa$. Now assuming $\mbox{sat}(\I\upharpoonright A)\leqslant\kappa$ we shall define a subtree $\mathcal{T}$ of $2^{\leqslant\kappa}$ and a mapping $\varphi:\mathcal{T}\mapsto P(A)$. Let $\emptyset\in\mathcal{T}$, $\varphi(\emptyset)=A$ and assume that levels $<\alpha$ are already constructed. If $\alpha=\beta+1$ for some $\beta<\kappa$ then proceed as follows. If $\varphi(t)\in\Ip$ for some $t\in\mbox{Lev}_\beta(\mathcal{T})$ then put $t\,\hat{}\,0,t\,\hat{}\,1$ into $\mathcal{T}$ and define a partition $\{\varphi(t\,\hat{}\,0),\varphi(t\,\hat{}\,1)\}\subseteq\Ip$ of $\varphi(t)$ (it is possible since $\I$ is atomless); if however $\varphi(t)\in\I$ put only $t\,\hat{}\,0$ in $\mathcal{T}$ and define $\varphi(t\,\hat{}\,0)=\varphi(t)$. For limit $\alpha$ put $t\in\mbox{Lev}_\alpha(\mathcal{T})$ if $t\upharpoonright\beta\in\mathcal{T}$ for all $\beta<\alpha$ and for such $t$ put $\varphi(t)=\bigcap_{\beta<\alpha}\varphi(t\upharpoonright\beta)$. Notice that by the above construction as well as by the fact that $\I\upharpoonright A$ is $\kappa$-saturated (each branch of $\mathcal{T}$ have less than $\kappa$ many splitting nodes) we can assume that
$$\varphi\left[\mbox{Lev}_{\kappa}(\mathcal{T})\right]\subseteq\I\quad\mbox{and}\quad A=\bigcup\varphi\left[\mbox{Lev}_{\kappa}(\mathcal{T})\right].$$
Therefore $$\left|\varphi\left[\mbox{Lev}_{\kappa}(\mathcal{T})\right]\right|\geqslant\mbox{add}(\I)=\kappa.$$
Thus $\mathcal{T}$ has at least $\kappa$ many branches and so does a tree $\varphi(\mathcal{T})\cap(\I\upharpoonright A)^+$. Pick
$$\left\{b_\alpha\in\left[\varphi(\mathcal{T})\cap(\I\upharpoonright A)^+\right]:\alpha<\kappa\right\}.$$
Let $\mathcal{D}$ be any dense subset of $(\I\upharpoonright A)^+$. We are going to show that a function $f:\mathcal{D}\mapsto\kappa$ defined as
$f(D)=\min\left\{\alpha<\kappa:D\In C\,\,\mbox{for some}\,\,C\in b_\alpha\right\}$
is unbounded. Then $|\mathcal{D}|\geqslant cf(\kappa)=\kappa$, so $\pi(\I\upharpoonright A)\geqslant\kappa$ and the proof will be complete.
We argue by a contradiction. Suppose $f\left[\mathcal{D}\right]\subseteq\lambda$ for some $\lambda<\kappa$. Pick any $C\in b_\lambda\setminus\bigcup\{b_\alpha:\alpha<\lambda\}$ and find $D\in\mathcal{D}$ such that $D\In C$. Now, on the one hand, by the definition of $f$, we have $C\cap B\geqslant_{\tiny\mathcal{I}}D\in\Ip$ for some $B\in b_{f(D)}$. On the other hand, by the assumption $f(D)<\lambda$ and by the fact that any element of $b_\lambda\setminus\bigcup\{b_\alpha:\alpha<\lambda\}$ is disjoint from any element of $\bigcup\{b_\alpha:\alpha<\lambda\}$, we have $C\cap B\in\I$. $\Box$

\medskip

Let us make some additional remarks. We say that a Boolean algebra $\mathbbm{B}$ is $\kappa$-represented (by $\{\mathbbm{B}_\alpha:\alpha<\kappa\}$) if it is co-absolute with a sum of $\kappa$ many Boolean algebras $\mathbbm{B}_\alpha$ such that $\pi(\mathbbm{B}_\alpha)<\kappa$. It is expected that $\kappa$-represented Boolean algebras should have density $\kappa$. For example we can extract from the proof of {\bf\em Theorem 1} the following. If in $P(X)/\I$ there exists a Boolean algebra $\kappa$-represented by {\em almost orthogonal} Boolean algebras (like $\langle b_\alpha\rangle/\I$ in the proof of {\bf\em Theorem  1}), then its density is at least $\kappa$. On the other hand the following theorem gives necessary condition for $P(X)/\I$ to be $\kappa$-represented by a chain of Boolean algebra.

\medskip

\indent {\bf\em Theorem 2} Let $\I$ be an atomless ideal on $X$ and let $\mathbbm{B}$ has a dense subalgebra $\mbox{add}(\I)$-represented by an increasing chain of regular subalgebras $\{\mathbbm{B}_\alpha:\alpha<\mbox{add}(\I)\}$. Then an inequality $\mbox{sat}(\I)\leqslant\mbox{add}(\I)$ forbids $P(X)/\I$ and $\mathbbm{B}$ to be equivalent.

\newpage

\noindent  {\bf\em Proof:}

Denote $\kappa:=\mbox{add}(\I)$ and suppose that $\mathbbm{B}$ is isomorphic to $P(X)/\I$. From the assumption of $\kappa$-representability and from {\bf\em Theorem 1} we conclude that for $\alpha<\kappa$
$$(\ast)\quad\pi\left(\mathbbm{B}_\alpha\right)<\kappa\leqslant\pi(\mathbbm{B}).$$
Using $(\ast)$ we shall construct a subtree $\mathcal{T}$ of a tree $2^{\leqslant\kappa}$ along with a family $\{a_t\in\mathbbm{B}^+:t\in\mathcal{T}\}$ being a Cantor (sub)scheme that satisfy the following
$$\quad(i)_\kappa\,\,\,\{a_t\in\mathbbm{B}^+:t\in\mbox{Lev$_\alpha(\mathcal{T})$}\}\,\,\mbox{is a partition in}\,\,\mathbbm{B}_\alpha\,\,\mbox{for all}\,\,\alpha<\kappa;$$
$$\quad(ii)_\kappa\,\,\,|\mbox{Split(Lev$_\alpha(\mathcal{T})$)}|\leqslant1\,\,\mbox{for all}\,\,\alpha<\kappa;\hspace{4.5cm}$$
$$\quad(iii)_\kappa\,\,\,|\mbox{Split($\mathcal{T}$)}|=\kappa.\hspace{8.2cm}$$
Here we use Split($\mathcal{A}$) to denote a family of splitting nodes in a subset $\mathcal{A}\subseteq\mathcal{T}$.

Let us argue that from the properties above it follows that there exists an antichain in $\mathbbm{B}$ of size $\kappa$. Indeed if there exists a branch $\mathcal{B}\in\left[\mathcal{T}\right]$ with $\mbox{Split($\mathcal{B}$)}=\{t_\alpha\in\mathcal{B}:\alpha<\kappa\}$ of size $\kappa$ then a family $\{b_{t_\alpha}-b_{t_{\alpha+1}}\in\mathcal{B}^+:\alpha<\kappa\}$ witnesses $\mbox{sat}(\I)=\mbox{sat}(\mathbbm{B})>\kappa=\mbox{add}(\I)$. If however each branch of $\mathcal{T}$ does not possess $\kappa$ many splitting nodes then $\mbox{ht}(\mathcal{T})\leqslant\kappa$ and we can extract from $\mbox{Split($\mathcal{T}$)}$ an antichain of size $\kappa$ as follows. Let $\alpha<\kappa$ and suppose we have already defined a partial antichain $\{t_\beta\in\mathcal{T}:\beta<\alpha\}$ in such a way that $\mbox{Split($\{s\in\mathcal{T}:t_\beta<s\}$)}=\emptyset$ and $\mbox{ht}(t_\beta)<\mbox{ht}(t_\gamma)<\mbox{ht}(\mathcal{T})\leqslant\kappa$ for all $\beta<\gamma<\alpha$. Consider a family $$\mathcal{S}_\alpha:=\left\{s\in\mbox{Split($\mathcal{T}$)}:\mbox{ht}(s)<\sup\{\mbox{ht}(t_\beta)<\kappa:\beta<\alpha\}+1\right\}.$$
Observe that $(ii)_\kappa$ and $(iii)_\kappa$ imply that $|\mbox{Split($\mathcal{T}$)}\setminus\mathcal{S}_\alpha|=\kappa$. When the last fact is combined again with $(iii)_\kappa$ and the assumption that there is no branch with $\kappa$ many splitting nodes then find the least $\gamma_\alpha<\kappa$ such that $\mbox{Lev$_{\gamma_\alpha}(\mathcal{T})$}\cap\mbox{Split($\mathcal{T}$)}\setminus\mathcal{S}_\alpha=\{s_{\gamma_\alpha}\}$ and define
$$t_\alpha=\sup\{t\in\mbox{Split($\mathcal{T}$)}\setminus\mathcal{S}_\alpha:s_{\gamma_\alpha}\leqslant t\}.$$
Observe that by the regularity of $\kappa$ a height of $t_\alpha$ is less than $\kappa$, so that it is a well-defined element of $\mathcal{T}$. Now, it is easily seen that this finishes an induction step. Hence $\{t_\alpha\in\mathcal{T}:\alpha<\kappa\}$ is an antichain in $\mathcal{T}$ of size $\kappa$ and by being a Cantor (sub)scheme a family $\{a_{t_\alpha}\in\mathbbm{B}^+:\alpha<\kappa\}$ witnesses $\mbox{sat}(\I)>\kappa$.

The tree $\mathcal{T}$ and the family $\{a_t\in\mathbbm{B}^+:t\in\mathcal{T}\}$ is constructed by an induction on levels of $\mathcal{T}$. Begin with $\mathcal{T}\upharpoonright 1=\{\emptyset\}$ and $a_\emptyset=\mathds{1}\in\mathbbm{B}_0$. Let $\alpha<\kappa$ and suppose we have already built $\mathcal{T}\upharpoonright\alpha$ and $\{a_t\in\mathbbm{B}^+:t\in\mathcal{T}\upharpoonright\alpha\}$ satisfying $(i)_\alpha$ and $(ii)_\alpha$.
Consider a family $$\mathcal{L}_\alpha:=\left\{\mathcal{B}\in\left[\mathcal{T}\upharpoonright\alpha\right]:\prod_{t\in\mathcal{B}}a_t\neq\mathds{O}\right\}.$$
First, observe that $\mathcal{L}_\alpha$ makes sense meaning that $\prod$ requires no subscripts since $\{a_t\in\mathbb{B}^+:t\in\mathcal{T}\upharpoonright\alpha\}\subseteq\mathbb{B}_\alpha$ and $\mathbb{B}_\alpha$ is a regular subalgebra of $\mathbb{B}$. Next, $$\mathcal{A}_\alpha:=\left\{\prod_{t\in\mathcal{B}}a_t\in\mathbb{B}_\alpha^+:\mathcal{B}\in\mathcal{L}_\alpha\right\}$$
is a pairwise disjoint family thanks to $(i)_\alpha$. Moreover, it follows from $(ii)_\alpha$ that $|\mathcal{T}\upharpoonright\alpha|\leqslant|\alpha|<\kappa$. Thus $\mathcal{A}_\alpha$ is a partition in $\mathbb{B}_\alpha$ by $\kappa$-completeness of $\mathbb{B}$. By $(\ast)$ and the fact that $\bigcup\{\mathbb{B}_\alpha^+:\alpha<\kappa\}$ is dense in $\mathbb{B}$ it follows that there exists $\beta\in(\alpha,\kappa)$ and $b\in\mathbb{B}_\beta$ such that $a_t-b\neq\mathds{O}$ for all $t\in\mathcal{T}\upharpoonright\alpha$. Let $\overline{\alpha}$ be the least such $\beta$ and fix $b\in\mathbb{B}_{\overline{\alpha}}$ as above. Since $\mathcal{A}_\alpha$ is a partition choose $a\in\mathcal{A}_\alpha$ meeting $b$, i.e. $a\cdot b\neq\mathds{O}$. Having all of this define $\mathcal{T}\upharpoonright[\alpha,\overline{\alpha}+1]$ and $\{a_t\in\mathbbm{B}^+:t\in\mathcal{T}\upharpoonright[\alpha,\overline{\alpha}+1]\}$ as follows $$\mathcal{T}\upharpoonright[\alpha,\overline{\alpha}]=\left\{\sup\mathcal{B}\,\hat{}\,(0)^\beta:\alpha\leqslant\sup\mathcal{B}+\beta\leqslant\overline{\alpha},\,\mathcal{B}\in\mathcal{L}_\alpha\right\},$$
$$a_t=\prod_{s\in\mathcal{B}}a_s\,\,\mbox{for}\,\,t\in\mathcal{T}\upharpoonright[\alpha,\overline{\alpha}]\,\,\mbox{with}\,\,t\upharpoonright\sup\mathcal{B}=\mathcal{B}\,\,\mbox{where}\,\,\mathcal{B}\in\mathcal{L}_\alpha.$$
and for $\overline{t}\in\mbox{Lev}_{\overline{\alpha}}(\mathcal{T})$ such that $a_{\overline{t}}=a$ put
$$\mbox{Lev}_{\overline{\alpha}+1}(\mathcal{T}):=\left\{t\,\hat{}\,(0)\in 2^{\overline{\alpha}+1}:t\in\mbox{Lev}_{\overline{\alpha}}(\mathcal{T})\setminus\{\overline{t}\}\right\}\cup\{t\,\hat{}\,(0),t\,\hat{}\,(1)\},$$
$$a_{t\,\hat{}\,(0)}=t\,\,\mbox{for}\,\,t\in\mbox{Lev}_{\overline{\alpha}}(\mathcal{T})\setminus\{\overline{t}\}\,\,\mbox{and}\,\,a_{\overline{t}\,\hat{}\,(0)}=a\cdot b, \,\,a_{\overline{t}\,\hat{}\,(1)}=a-b.$$
and thus $\overline{t}$ is a splitting node. This finishes the induction step. It is now immediate that $(i)_\kappa$ and $(ii)_\kappa$ are satisfied. Moreover since by construction for any $\alpha<\kappa$ we can find a splitting node on level $\overline{\alpha}>\alpha$ the last property $(iii)_\kappa$ also holds. The proof of theorem is completed. $\quad\Box$

\medskip

Recall that a Boolean algebra $C(\lambda)$ is called {\em Cohen algebra} (\cite{J}) if it is equivalent to the forcing adding $\lambda$ many Cohen reals, i.e. it includes as a dense subset a poset of finite binary functions with domain being subset of $\lambda$ and ordering being a reverse inclusion. Thus by the definition $\pi(C(\lambda))=\lambda$. Moreover, after obvious identification of $C(\lambda')$ with a regular subalgebra of $C(\lambda)$ for $\lambda'<\lambda$, an algebra $\bigcup\{C(\lambda'):\lambda'<\lambda\}$ is co-absolute with $C(\lambda)$. Therefore by {\bf\em Theorem 2} w can conclude

\indent {\bf\em Corollary 3} If $P(X)/\I$ is equivalent to $C(\lambda)$ then $\lambda>\mbox{add}(\I)$.

\medskip

Similarly (using almost the same idea) one can prove analogical result for measure algebras or even more general Boolean algebras. Before doing this let us establish the following simple lemmas.

\medskip

\indent {\bf\em Lemma 4} If $\I$ is an ideal with $\mbox{sat}(\I)\leqslant\mbox{add}(\I)$ then $wh(\I)\geqslant\mbox{add}(\I)$.

\medskip

\indent  {\bf\em Proof:}

Let $\alpha<\mbox{add}(\I)$ and let $\{\mathcal{P}_\beta:\beta<\alpha\}$ be an arbitrary matrix of partition in $P(X)/\I$ with enumerations $\mathcal{P}_\beta=\{[A_{\beta\xi}]:\xi<\gamma_\beta\}$. Note that $\gamma_\beta<\mbox{add}(\I)$ by the assumption. Define the following pairwise disjoint family $$\mathcal{P}:=\left\{\prod_{\beta<\alpha}[A_{\beta f(\beta)}]\in\left(P(X)/\I\right)^+:f\in\prod_{\beta<\alpha}\gamma_\beta\right\}.$$
Observe that $\mathcal{P}$ refines $\{\mathcal{P}_\beta:\beta<\alpha\}$ and moreover it is a partition as
$$\sum\mathcal{P}=\sum\left\{\left[\bigcap_{\beta<\alpha}A_{\beta f(\beta)}\right]:f\in\prod_{\beta<\alpha}\gamma_\beta\right\}=$$
$$=\left[\bigcup\left\{\bigcap_{\beta<\alpha}A_{\beta f(\beta)}:f\in\prod_{\beta<\alpha}\gamma_\beta\right\}\right]=\left[\bigcap_{\beta<\alpha} \bigcup_{\xi<\gamma_\beta}A_{\beta\xi}\right]=$$
$$=\prod_{\beta<\alpha}\left[\bigcup_{\xi<\gamma_\beta}A_{\beta\xi}\right]=\prod_{\beta<\alpha}\left(\sum\mathcal{P}_\beta\right)=\mathds{1},$$
where in the first and in the third line we use $\max\{\alpha,\gamma_\beta,|\mathcal{P}|\}<\mbox{add}(\I)$ together with $\mbox{add}(\I)$-completeness of $P(X)/\I$ and in the second line we use complete distributivity of $P(X)$. Finally, as the partition $\mathcal{P}$ refines all partitions from $\{\mathcal{P}_\beta:\beta<\alpha\}$ we obtain $\alpha<wh(\I)$. Hence $wh(\I)\geqslant\mbox{add}(\I)$ since $\alpha$ was arbitrary. $\Box$

\medskip

\indent {\bf\em Lemma 5} If $\I$ is an atomless ideal on $X$ then $wh(\I)\leqslant\pi(\I)$.

\medskip

\indent  {\bf\em Proof:}

\medskip

First, let $\{D_\alpha\in\I^+:\alpha<\pi(\I)\}$ be dense in $\I^+$ and define a mapping $\varphi:\Ip\mapsto\pi(\I)$ by the rule
$\varphi(A)=\min\left\{\alpha<\pi(\I):D_\alpha\subseteq_\I A\right\}.$ Observe that
\begin{itemize}
  \item if $A_1\supseteq_\I A_0\in\Ip$ then $\varphi(A_1)\leqslant\varphi(A_0)$;
  \item for all $A\in\Ip$ there is $A\supseteq_\I B\in\Ip$ such that $\varphi(A)<\varphi(B)$.
\end{itemize}
The first observation is a direct consequence of the definition of $\varphi$; for the second, just take for $B$ a set $E$ such that $\left\{E,D_{\varphi(A)}\setminus E\right\}$ is a partition of $D_{\varphi(A)}$ (the choice is possible since $\I$ is atomless).

Suppose that $\pi(\I)<wh(\I)$ and proceed as follows to get a contradiction. We shall construct a matrix of partitions $\{\mathcal{P}_\alpha\subseteq\Ip:\alpha<\pi(\I)\}$ with the following properties
%\,\,\mathcal{P}_0=\{B_0,X\setminus B_0\}\cap\Ip\,\,
%$$\quad(i)\,\,\,\mbox{for all}\,\,\beta<\alpha\,\,\mbox{there is}\,\,A\in\mathcal{P}_\beta\,\,\mbox{with}\,\,A\cap B_\beta\in\Ip;$$
$$\quad(i)_{<\pi(\I)}\,\,\,\mbox{if}\,\,\beta<\gamma<\alpha\,\,\mbox{then}\,\,\mathcal{P}_\gamma\,\,\mbox{almost refines}\,\,\mathcal{P}_\beta\hspace{4.3cm},$$
$\hspace{2.1cm}\mbox{i.e. if}\,\,A\in\mathcal{P}_\gamma\,\,\mbox{then}\,\,A\cap C\in\Ip\,\,\mbox{for finitely many}\,\,C\in\mathcal{P}_\beta;$
$$\quad(ii)_{<\pi(\I)}\,\,\,\mbox{if}\,\,\beta<\gamma<\alpha,\,\,A\in\mathcal{P}_\beta,\,\,C\in\mathcal{P}_\gamma,\,\,A\cap C\in\Ip\,\,\mbox{then}\,\,\varphi(A)<\varphi(C).$$ Put $\mathcal{P}_0=\{X\}$, let $\alpha<\pi(\I)$ and suppose we have already built a partial matrix $\{\mathcal{P}_\beta\subseteq\Ip:\beta<\alpha\}$ satisfying $(i)_{<\alpha}$ and $(ii)_{<\alpha}$. We find appropriate $\mathcal{P}_\alpha$ as follows. First, choose a matrix $\mathcal{P}'_\alpha$ almost refining $\{\mathcal{P}_\beta:\beta<\alpha\}$; this is possible by $\alpha<wh(\I)$. Next refine each member $A\in\mathcal{P}'_\alpha$ to a partition $\mathcal{P}'_\alpha(A)$ of $A$ in the way that for any $B\in\mathcal{P}'_\alpha(A)$ it holds $\varphi(B)>\varphi(A)$; this is possible by the aforementioned  properties of $\varphi$ and by Kuratowski-Zorn Lemma. Now it is easily seen that $\mathcal{P}_\alpha:=\bigcup\{\mathcal{P}'_\alpha(A):A\in\mathcal{P}'_\alpha\}$ is a partition satisfying $(i)_\alpha$ and $(ii)_\alpha$ as demanded. Note also that $(ii)_{<\pi(\I)}$ implies that if $A\in\mathcal{P}_\alpha$ then $\varphi(A)\geqslant\alpha$.

As $\pi(\I)<wh(\I)$ let  $\mathcal{P}$ be a partition almost refining all $\{\mathcal{P}_\alpha:\alpha<\pi(\I)\}$. In particular there exists $E\in\mathcal{P}$ such that for each $\alpha<\pi(\I)$ there is $A_\alpha\in\mathcal{P}_\alpha$ such that $A_\alpha\cap E\in\Ip$. Take $\alpha\geqslant\varphi(E)$ and obtain (again using $(ii)_{<\pi(\I)}$) the following
$$\varphi(E)>\varphi(A_\alpha)\geqslant\alpha=\varphi(E).$$
The above contradiction completes the proof.  $\Box$

\medskip
 
As an easy consequence of the above lemmas we obtain

\medskip

\indent {\bf\em Theorem 6}

\noindent If $\I$ is an atomless ideal satisfying $\mbox{sat}(\I)\leqslant\mbox{add}(\I)$ and $h(\I)<wh(\I)$ then $\mbox{add}(\I)<\pi(\I)$.

\medskip

\indent  {\bf\em Proof:}

\medskip 

Straightforward, by $\mbox{sat}(\I)\leqslant\mbox{add}(\I)$ we have $h(\I)\geqslant\mbox{add}(\I)$ due to {\bf\em Lemma 4}. This together with the second assumption and {\bf\em Lemma 5} yields $$\mbox{add}(\I)\leqslant h(\I)<wh(\I)\leqslant\pi(\I).\quad\Box$$

\medskip

Recall that a Boolean algebra $\mathbbm{B}$ is called {\em measure algebra} if it carries a measure, i.e. non-decreasing, real-valued, non-negative, $\sigma$-additive, normalized ($\mu(\mathbbm{1})=1$) function $\mu$ vanishing only at $\mathds{O}$. A measure $\mu$ is atomless iff $\I_\mu:=\{a\in\mathbbm{B}:\mu(a)=0\}$ is an atomless ideal. It is well-known that measure algebra is non-distributive, weakly $\omega$-distributive Boolean algebra such that $\mbox{sat}(\I)=\omega_1$ and $h(\I)<wh(\I)$. A random algebra $R(\lambda)$ is a Boolean algebra co-absolute-with forcing adding $\lambda$ many random reals, that is, it includes as a dense subset a partial order $(B_\lambda,\subseteq)$, where $B_\lambda$ is a family of Borel subsets of $\{0,1\}^\lambda$ which have positive Haar measure. Since by definition $\pi(R(\lambda))\geqslant\lambda$, {\bf\em Theorem 6} implies

\medskip

\indent {\bf\em Corollary 7} If $P(X)/\I$ is equivalent to forcing with $\lambda$ many random reals then $\lambda>\mbox{add}(\I)$.

\medskip

\indent We finish with some remarks. The use of the language of trees is not circumstancial. This is perfectly visible when comparing tree forcing notions with forcing with ideal. Boolean algebras generated by these forcings are in a sense canonical, they are the most representative and interesting among all definable forcings \cite{Z}. Still, the most intriguing problem remains: what suitable definable forcings (Boolean algebras) are not equivalent with $P(X)/\I$? Some results appeared in \cite{GS} for the simplest tree forcings however the problem is widely open. In forthcoming paper we use the technique form the paper to systematic study so called idealized forcings \cite{Z} and its properties preventing to be equivalent to $P(X)/\I$.

 {\tiny\sc
$\begin{array}{cc}
 Ryszard Frankiewicz &  S{\l}awomir Szczepaniak \\
 Institute of Mathematics, Polish Academy of Sciences & Institute of Mathematics, Polish Academy of Sciences \\
  \'Sniadeckich 8, 00-950 Warszawa, Poland & \'Sniadeckich 8, 00-950 Warszawa, Poland \\
  e-mail: rfrankiewicz@mac.com &  e-mail: szczepaniak@impan.pan.wroc.pl
\end{array}$
}

%{\tiny\sc Ryszard Frankiewicz}\\ Institute of Mathematics, Polish
%Academy of Sciences,\\ \'Sniadeckich 8, 00-950 Warszawa, Poland \\
%{\em e-mail: rf@impan.pl}\\
%
%{\tiny\sc S{\l}awomir Szczepaniak}\\ Institute of Mathematics, Polish
%Academy of Sciences,\\ \'Sniadeckich 8, 00-950 Warszawa, Poland \\
%{\em e-mail: szczepaniak@impan.pan.wroc.pl}

\end{document}